\documentstyle[12pt,amsfonts]{article}
\newlength{\abstractwidth}
\renewcommand{\thefootnote}{\fnsymbol{footnote}}
\renewcommand{\thanks}[1]{\footnote{#1}} 
\newcommand{\starttext}{ \setcounter{footnote}{0}
\renewcommand{\thefootnote}{\arabic{footnote}}}

\newcommand{\be}{\begin{equation}}
\newcommand{\bea}{\begin{eqnarray}}
\newcommand{\eea}{\end{eqnarray}} \newcommand{\ee}{\end{equation}}
 \newcommand{\<}{\langle}
\renewcommand{\>}{\rangle} \def\ba{\begin{eqnarray}}
\def\ea{\end{eqnarray}}

\def\o{\omega}

\def\tr{{\rm tr}}
\def\det{{\rm det}}

\def\log{\,{\rm log}\,}

\def\o{\omega}

\def\o{\omega}

\def\na{\nabla}

\def\ge{\geq}
\def\le{\leq}

\def\qed{$\Box$}

\def\p{\partial}

\def\ddb{{\partial\bar\partial}}

\def\na{{\nabla}}

\def\[{{\bf [}}
\def\]{{\bf ]}}

\def\ddbar{\ddb}

\begin{document}
\starttext \baselineskip=18pt \setcounter{footnote}{0}
\newtheorem{theorem}{Theorem}
\newtheorem{lemma}{Lemma}
\newtheorem{corollary}{Corollary}
\newtheorem{definition}{Definition}
\newtheorem{conjecture}{Conjecture}
\newtheorem{proposition}{Proposition}

\begin{center}
{\Large \bf ON CONVERGENCE CRITERIA FOR THE COUPLED FLOW OF LI-YUAN-ZHANG
\footnote{Work supported in part by the National Science Foundation under grants DMS-12-66033 and DMS-17-10500.}}
\nonumber \\
\medskip
\centerline{
Teng Fei, Bin Guo, and Duong H. Phong}

\bigskip

\begin{abstract}

\medskip

{\small A one-parameter family of coupled flows
depending on a parameter $\kappa>0$ is introduced which reduces when $\kappa=1$ to the coupled flow of a metric $\o$ with a $(1,1)$-form $\alpha$ due recently to Y. Li, Y. Yuan, and Y. Zhang. It is shown in particular that, for $\kappa\not=1$, estimates for derivatives of all orders would follow from $C^0$ estimates for $\o$ and $\alpha$. Together with the monotonicity of suitably adapted energy functionals, this can be applied to establish the convergence of the flow in some situations, including on Riemann surfaces. Very little is known as yet about the monotonicity and convergence of flows in presence of couplings, and conditions such as $\kappa\not=1$ seem new and may be useful in the future.}

\end{abstract}

\end{center}

\baselineskip=15pt

\section{Introduction}
\setcounter{equation}{0}

The main goal of this paper is to study the following flow of a K\"ahler metric $\o_t$ coupled to a closed $(1,1)$-form $\alpha_t$, on a given compact complex manifold $X$,
\bea
\label{LYZo}
\p_t\o_t&=&-Ric(\o_t)+\lambda\o_t+\alpha_t,\\
\label{LYZa-theta}
\p_t\alpha_t&=& -\kappa\,\bar\p\bar\p^\dagger_{\o_t} \alpha_t.
\eea
Here $Ric(\o)=-i\p\bar\p\log \o^n$ is the Ricci form of $\o$, $n$ is the complex dimension of $X$, and $\kappa$ a positive constant. When $\kappa=1$ and $\lambda=-1$, this is the flow introduced recently by Li, Yuan, and Zhang \cite{LYZ}. However, the stationary points are independent of $\kappa$, and as we shall see in the present paper, it may be advantageous to consider a different value of $\kappa$. We shall refer to the flow (\ref{LYZo},\ref{LYZa-theta}) as the $\kappa$-LYZ flow. The original motivation for the $\kappa$-LYZ flow is that its stationary points are given by a K\"ahler metric of constant scalar curvature together with a harmonic $(1,1)$-form. As such, it provides a natural approach to the well-known open problem of finding K\"ahler metrics of constant scalar curvature \cite{LYZ}.

\medskip
Our interest in the $\kappa$-LYZ flow comes from a different source. A leading contender for a unified theory of the forces of nature at their most fundamental level is M theory and its limiting string and supergravity theories \cite{BBS, HW, T}. These theories all incorporate the gravitational field together with other fields, and the resulting equations, whether in the original dimension or upon compactification, are all Einstein's equation modified by interactions with other fields or higher string modes. The corresponding parabolic flows are then all essentially Ricci flows, modified by couplings to other fields. Some explicit examples are the renormalization group flow for the bosonic string considered in \cite{Li1}, the Anomaly flow arising from compactifications of the heterotic string \cite{PPZ1, PPZ2, FHP} and the flows found in \cite{FGP} arising from compactifications of eleven-dimensional supergravity.
However, while there has been considerable progress in recent years in the understanding of the Ricci flow, there are still very few tools available for the study of non-linear coupled systems in general, and in particular of the
long-time behavior of the Ricci flow when it is coupled to other fields. Even for
the simplest coupled Ricci flow, namely the Ricci flow coupled to a scalar field, only a criterion for the development of singularities \cite{List} and a Perelman-type pseudo-locality theorem \cite{GHP} are known. A generalization of the first of these results to the Ricci flow coupled to the harmonic map flow is in \cite{Mu}. The $\kappa$-LYZ flow is arguably the simplest example of a coupled Ricci flow in K\"ahler geometry, since the additional field $\alpha$ is a closed $(1,1)$-form, and the closedness of both the metric form $\o$ and of the additional field $\alpha$ are preserved along the flow. It can also be viewed as an Abelian model for the Anomaly flow \cite{PPZ1,PPZ2}, in which the role of $\alpha$ is played by a Hermitian metric on a holomorphic vector bundle, and higher powers of the curvature appear in the couplings.

\medskip

The study of a partial differential equation usually begins with the identification of a minimum number of estimates from which the existence and/or $C^\infty$ regularity of the equation would follow. Fundamental examples are the $C^{2,\alpha}$ estimates for linear uniformly elliptic second order equations, the $C^0$ estimate for the complex Monge-Amp\`ere equation on a compact K\"ahler manifold, and the uniform estimate for the metrics in the method of continuity for the problem of constant scalar curvature K\"ahler metric. In all these cases
estimates for the derivatives of any order would follow, implying in turn both existence and regularity properties for the equation. Even such basic results are not yet available for either general non-linear elliptic or parabolic systems, or even more specifically for the Ricci flow coupled with other fields.

\medskip
The main goal of the present paper is to initiate such a study for the $\kappa$-LYZ flow. One of our main results is that, for $\kappa\not=1$, if the metrics $\o_t$ are all equivalent and the form $\alpha_t$ are uniformly bounded, then the derivatives of both $\o_t$ and $\alpha_t$ of all orders are uniformly bounded (Theorem 1).
Thus the uniform boundedness of $\o_t$ and $\alpha_t$ plays the role in the $\kappa$-LYZ system of the $C^{2,\alpha}$ estimate for linear elliptic equations, the $C^0$ estimate for the Monge-Amp\`ere equation, and the uniform estimates for metrics in the continuity method for the problem of constant scalar curvature.
The condition $\kappa\not=1$ seems to be new and it is essential to our proof. It can be viewed as a restriction on how the Ricci flow can be coupled to the flow of the form $\alpha_t$, and such restrictions may be important in the future investigation of other coupled systems.
The uniform bounds on $\o_t$ and $\alpha_t$ in Theorem 1 are strong conditions, but we can show that they hold for suitable data and suitable values of $\kappa$ in the case of Riemann surfaces (Theorem 2). A key tool in this case is an apparently new Liouville-Entropy type energy for the coupled system which is monotone along the $\kappa$-LYZ flow. We can then establish the convergence of the $\kappa$-LYZ flow in this case. While this result is rather special, convergence results are rather rare in general for coupled systems, and this is to our knowledge the only result on the convergence of a coupled Ricci flow available so far. When the manifold $X$ has higher dimension, we introduce instead a modified Mabuchi energy which incorporates the additional field, and which can be shown to be monotone along the $\kappa$-LYZ flow. We can then show that, under the same hypotheses as in Theorem 1 and the additional assumption that the modified Mabuchi energy is bounded from below, the $\kappa$-LYZ flow converges smoothly to a smooth stationary point (Theorem 3).

\medskip
It may be instructive to relate our results to some recent ones in the literature. In the original paper of Li, Yuan, and Zhang \cite{LYZ}, Shi-type estimates were established for all derivatives of the curvature, assuming that the curvatue is bounded. But Shi-type estimates hold only for finite-time intervals, while we are here interested in bounds that are uniform for all time. Another closely related problem is the problem of K\"ahler metrics of constant scalar curvature, which has already been mentioned several times. It is an elliptic equation of 4th-order in the potential, which can also be expressed as a coupled elliptic system in the potential and the volume, viewed as two separate unknowns\footnote{A similar decomposition was used by Trudinger and Wang \cite{TW} for the affine Plateau problem.}. It is then well-known that the condition of uniform boundedness for the metrics in the method of continuity would imply estimates of all orders and the solvability of the equation. The extension of this result to the parabolic case does not seem available in the literature as yet. Additional difficulties result from the need to estimate the partial derivatives in time, and the resulting complicated mixing of the two unknowns. Theorem 1 in the present paper can be interpreted as such a parabolic result for the $\kappa$-LYZ equation. Very recently,
Chen and Cheng \cite{CC1, CC2, CC3} made a breakthrough on the constant scalar curvature problem by showing in particular how the uniform boundedness of the metrics in the elliptic approach follows from the properness of the Mabuchi functional. It is however not known at this moment how to adapt their methods to the parabolic setting. For example, the classical parabolic approach to the constant scalar curvature problem, namely by the Calabi flow, still remains an open problem at this time.

\section{Estimates for Higher Order Derivatives of $\o$ and $\alpha$}
\setcounter{equation}{0}

We begin with some simplifying remarks on the normalization of the $\kappa$-flow (\ref{LYZo},\ref{LYZa-theta}). Clearly, if $(\o_\infty,\alpha_\infty)$ is a stationary point of the flow, then the K\"ahler classes $[\o_\infty],[\alpha_\infty]$ must satisfy
\bea
-2\pi c_1(X)+\lambda[\o_\infty]+[\alpha_\infty]=0.
\eea
It is natural then to consider initial data $\o_0,\alpha_0$ in cohomology classes satisfying the same condition. Since the flow preserves the cohomology class of $\alpha_0$, it means that we should assume the following condition on the initial data
\bea
\label{nec}
-2\pi c_1(X)+\lambda[\o_0]+[\alpha_0]=0.
\eea
In general, it is not difficult to work out explicitly the dependence of $[\o_t]$ on time, and to see that (\ref{nec}) is also a necessary condition for the convergence of $[\o_t]$ when $\lambda\geq 0$. Thus it is simplest to assume henceforth the condition (\ref{nec}) on the initial data $(\o_0,\alpha_0)$. As a consequence $[\o_t]$ is constant in time. Moreover, the $\kappa$-LYZ flow is invariant under the scaling $\omega\to M\omega$, $t\to Mt$, $\alpha\to\alpha$, $\lambda\to \lambda/M$, therefore we may always assume that
\bea
V=\int_X{\o_t^n\over n!}=1.
\eea
In this section, we shall prove the following theorem:

\begin{theorem}
\label{higherorder}
Assume that the $\kappa$-LYZ flow exists on a time interval $[0,T)$ and that
the following basic assumption holds: there exists a positive constant $K > 1$ so that
\bea\label{eqn:assumption}
K^{-1} \omega_0 \le \omega_t \le K \omega_0,\quad - K \omega_t \le \alpha_t \le K \omega_t
\eea
for all $t\in [0,T)$. Let the endomorphisms $h_t$ be defined by $h^p{}_q=\o^{p\bar m}_0\o_{\bar m q}$. If $0 < \kappa \neq 1$, then for any integer $k$, there exists a constant $C_k$ depending only on $\kappa,K,k,n,\lambda$ and the initial data $\o_0$ and $\alpha_0$ so that
\bea
|\na^k h|+|\na^k \alpha|\leq C_k.
\eea
\end{theorem}

\medskip

These estimates are similar in spirit to the ones for K\"ahler-Einstein metrics, in the sense that all estimates in that case follow from the $C^0$ estimate. The difference resides in the fact that the estimates in the K\"ahler-Einstein case are estimates for the potential, while here they are for the metrics and $(1,1)$-form $\alpha$. The fact that we are dealing here with a system also creates many new difficulties.

\medskip
As in \cite{LYZ}, we shall make use of the reformulation of the flow in terms of potentials. Thus set
\bea\label{eqn:potential}\omega_t = \omega_0 + i\ddb \varphi_t, \quad \alpha_t = \alpha_0 - i\ddb F_t,\eea
for smooth functions $\varphi_t$ and $F_t$ determined up to additive terms which are constants in space, but which may depend on time. We can fix these additive terms by requiring that the $\kappa$-LYZ be equivalent to the following coupled flow of the system $(\varphi_t,F_t)$,
\bea
\label{eqn:varphi}
\frac{\partial \varphi}{\partial t} & =& \log\Big( \frac{\omega_t^n}{\omega_0^n} \Big)- H_0 + \lambda \varphi - F,\\
\frac{\partial F}{\partial t} & = & \kappa\,\Delta F - \kappa\,\tr_{\omega_t} \alpha_0 + \kappa\, b.
\label{eqn:F}
\eea
Here $H_0\in C^\infty$ is the Ricci potential defined by
\bea
Ric(\omega_0) - \lambda \omega_0 -\alpha_0 = i\ddb H_0
\qquad
\int_X e^{H_0} \omega_0^n = \int_X \omega_0^n = V,
\eea
and $b$ is a time-independent constant chosen so that
\bea
\label{b}
\int_X \dot F \omega_t^n = 0,\quad {\mathrm{hence}} ~ b = n  \int_X \alpha_0 \wedge \omega_0^{n-1}.
\eea
In this form, the $\kappa$-LYZ flow with $\lambda=0$ and $\alpha$ the Ricci form of some metric is a parabolic version of the elliptic system of equations for the potentials of $\o$ and $\alpha$ considered by Chen and Cheng \cite{CC1, CC2, CC3}.

\medskip

By the assumption (\ref{eqn:assumption}), we know that $\ddb \varphi_t$ and $\ddb F_t$ are both uniformly bounded with respect to the fixed K\"ahler metric $\omega_0$.  Our goal is to derive higher order derivative estimates of $\omega_t$ and $\alpha_t$, beginning with the third order derivative estimate for $\varphi_t$. For the Monge-Amp\`ere equation, $C^3$ estimates for the potential are obtained using the maximum principle and the Calabi identity (see e.g. \cite{Y}). Here the Calabi identity and the maximum principle do not suffice, because the coupling between $\o_t$ and $\alpha_t$ results in terms that involve the derivatives of $\alpha_t$ and hence are not a priori bounded. We overcome these difficulties by combining the method of \cite{PSS}, which considers instead the evolution of the connection defined by $\o_t$, with a parabolic Moser iteration argument. For this parabolic Moser argument, it is crucial to have a uniform lower bound for the scalar curvature of $\o_t$, and this is where the condition $\kappa\not=1$ is needed.

\medskip

To begin with, we establish a lower bound for the scalar curvature $R$.

\begin{lemma}\label{lemma 1}
For any $0<\kappa\neq 1$, there exists a constant $C=C(\omega_0,\kappa,K,\lambda)>0$ such that $$R \ge -C,$$
where $R = R(\omega_t)$ is the scalar curvature of $\omega_t$.
\end{lemma}
\noindent{\it Proof.}
From the definition of the $\kappa$-LYZ flow, it is easy to derive the flows for $R$ and $\tr_\omega\alpha$,
\bea\nonumber
\frac{\partial R}{\partial t}& = & \frac{\partial }{\partial t} ( - g^{i\bar j}\partial_i \partial_{\bar j} \log \det g  )\\
& = & \nonumber - \Delta ( -R + n\lambda + \tr_\omega\alpha  ) + R_{\bar j i}(R_{\bar i j} - \lambda g_{\bar i j} - \alpha_{\bar i j})\\
& = & \label{eqn:for scalar}  \Delta R - \Delta \tr_\omega\alpha + |Ric|^2- \lambda R - R_{\bar j i} \alpha_{\bar i j}
\eea
while
\bea\nonumber
\frac{\partial}{\partial t} \tr_\omega\alpha & = &\kappa\,\Delta \tr_\omega\alpha + (R_{\bar i j} - \lambda g_{\bar i j} - \alpha_{\bar i j}) \alpha_{\bar j i}.
\eea
For $\kappa\neq 1$, we can combine these two equations to get
\bea
\label{kappa}
\nonumber
 \left(\frac{\partial}{\partial t} - \Delta \right)\left( R+ \frac{\tr_\omega\alpha}{\kappa -1} \right) &= & |Ric|^2- \lambda R - \lambda \tr_\omega\alpha - |\alpha|^2 - \frac{\kappa-2}{\kappa-1}R_{\bar j i}\alpha_{i\bar j}\\
 & \ge & \frac{R^2}{2n} - \lambda R - C\nonumber\\
 & \ge & \frac{1}{4n}\left(R+\frac{\tr_\omega\alpha}{\kappa-1}\right)^2  - C\nonumber,
\eea
where in the first inequality we use the fact that $n|Ric|^2\ge R^2$.
We can apply now the maximum principle to the function $R+\tr_\omega\alpha/(\kappa-1)$ and obtain a uniform lower bound for $R+\tr_\omega\alpha/(\kappa-1)$. Since by assumption $\tr_\omega\alpha$ is bounded, the lower bound for $R$ follows.

\medskip

\noindent We are now ready to prove the Calabi $C^3$-estimates of $\varphi$, under the assumption (\ref{eqn:assumption}).

\begin{lemma}\label{lemma 2}
There exists a constant $C=C(\omega_0, \kappa, K, \lambda)>0$ such that
$$\sup_{X\times [0,T)}|\nabla \ddb\varphi|_{\omega_0} \le C.$$
\end{lemma}
\noindent {\it Proof.}
For notation convenience we will denote $\hat g$ the metric associated to the fixed K\"ahler form $\omega_0$. We define
$$S^{k}_{ij}: = \Gamma^k_{ij} - \hat\Gamma^k_{ij}$$
where $\Gamma^k_{ij} = g^{k\bar l} \frac{\partial g_{\bar l i}}{\partial z^j}$ is the Christoffel symbol of $\omega = \omega_t$ and $\hat \Gamma^k_{ij}$ is that of $\hat g$. Because of the equivalence of $\omega_t$ and $\omega_0$ by the assumption (\ref{eqn:assumption}), to prove the lemma it suffices to show that $|S|^2_{\omega} = S^k_{ij} \overline{S^{r}_{pq}} g^{i\bar p} g^{j\bar q} g_{\bar r k}$ is uniformly bounded.

We calculate under the normal coordinates of $\omega = \omega_t$ at some fixed point,
\bea\nonumber
\frac{\partial}{\partial t} |S|^2 & = & 2 (R_{\bar i p} - \lambda g_{\bar ip} - \alpha_{\bar i p}) S_{ij}^k \overline{S^k_{pj}} + (- R_{\bar r k} + \lambda g_{\bar r k} + \alpha_{\bar r k}) S^k_{ij} \overline{S^r_{ij}} + 2 Re \Big(\frac{\partial S^k_{ij}}{\partial t} \overline{S^k_{ij}} \Big)\\
& = & 2 (R_{\bar i p} - \lambda g_{\bar ip} - \alpha_{\bar i p}) S_{ij}^k \overline{S^k_{pj}} + (- R_{\bar r k} + \lambda g_{\bar r k} + \alpha_{\bar r k}) S^k_{ij} \overline{S^r_{ij}}\nonumber\\
 && + 2 Re \Big(\nabla_j (- R_{\bar k i} + \alpha_{\bar k i})\overline{S^k_{ij}}  \Big).    \nonumber
\eea
And if we denote $S^k_{ij, p} = \nabla_p S^k_{ij}$ and $S^k_{ij,p\bar p} =\bar  \nabla _{\bar p } \nabla_p S^k_{ij}$, then
\bea\nonumber
\Delta |S|^2 & = & (S^k_{ij}\overline{S^k_{ij}})_{p\bar p}\\
& = & \nonumber S^k_{ij,p\bar p } \overline{S^k_{ij}} + S^k_{ij,p} \overline{ S^k_{ij,p}} +  S^k_{ij,\bar p}\overline{ S^k_{ij,\bar p}} + S^k_{ij} \overline{S^k_{ij,\bar p p }}\\
& = &\overline{S^k_{ij}} ( S^k_{ij,\bar p p} + S^k_{mj} R_{\bar m i} + S^k_{im} R_{\bar m j} - S^m_{ij} R_{\bar k m}) + S^k_{ij} \overline{S^k_{ij,\bar p p }  } + |\nabla S|^2 + |\bar \nabla S|^2. \nonumber
\eea
Note that
\bea\nonumber
S^k_{ij,\bar p p}& = & \nabla _p \bar \nabla_{\bar p} S^k_{ij} = \nabla_p ( \bar \partial_{\bar p} ( \Gamma^k_{ij} - \hat\Gamma ^k_{ij} )    )\\
& = & \nabla _p ( - R_{~i\bar p j}^k + \hat R^k_{~ i\bar p j}  )\nonumber\\
& = & - \nabla_j R^k_{~ i} + ( \nabla_p - \hat \nabla _p  ) \hat R ^k_{~i\bar p j} + \hat \nabla_j \hat R^k_{~ i}\nonumber
\eea
and the middle term on the RHS can be written as the form $S* \hat Rm$.

Therefore we have
\bea\nonumber
( \frac{\partial}{\partial t} - \Delta) |S|^2 & = & - 2 \alpha_{\bar i p} S^k_{ij} \overline{S^k_{pj}}  -  \lambda |S|^2 + \alpha_{\bar r k} S^k_{ij} \overline{S^r_{ij}}+ 2 Re ( \nabla_j \alpha_{\bar k i} \overline{S^k_{ij}}   ) \\
& & \nonumber  - 2 Re( (\nabla_p - \hat \nabla_p ) \hat R^k_{~i\bar p j} \overline{S^k_{ij}}   ) - 2 Re ( \hat \nabla _j \hat R^k_{~i}  \overline{S^k_{ij}}  ) - |\nabla S|^2 - |\bar \nabla S|^2\\
& \le & C |S|^2  + C |S| + 2 Re( \nabla_j \alpha_{\bar k i} \overline{S^k_{ij}}  ) - |\nabla S|^2 - |\bar \nabla S|^2,    \label{eqn:S sqr}
\eea for some $C =C(\omega_0, K, \lambda)>0$.
Since $\nabla_j \alpha_{\bar k i}$ is not apriorily bounded, we cannot apply maximum principle as usual to estimate $|S|^2$. Instead we will bound $|S|^2$ by Moser iteration. We denote $\hat u = \max( |S|^2, 1  )\ge 1$ and for any times $ t_0 < t_1  \le T' < T  $ define a Lipchitz function $\xi(t)$ such that $\xi(t) = 0$ for $t\le t_0$, $\xi(t) = 1$ for $t\ge t_1$ and $\xi(t) = \frac{t-t_0}{t_1 - t_0}$ for $t\in [t_0,t_1]$. Multiplying both side of (\ref{eqn:S sqr}) by $\xi(t) \hat u^p$ for $p\ge 1$, we get for any $s\in [t_1, T]$
\bea\nonumber
&& \int_X \frac{\xi \hat u^{p+1}}{p+1} \omega_t^n \Big|_{t=s} + \int_{t_1}^s \int_X \frac{4p }{(p+1)^2} |\nabla \hat u^{\frac{p+1}{2}}|^2 \omega_t^n dt\\
&\le & \nonumber \int_{t_0}^T \int_X\Big(C  \hat u^{p+1} + 2 \xi \hat u^p Re(\nabla_j \alpha_{\bar k i}\overline{S^k_{ij}}) - \xi \hat u^p(|\nabla S|^2 + |\bar \nabla S|^2)  \Big)\omega_t^n dt\\
&& \nonumber +\int_{t_0}^T \xi'(t) \int_X \hat u^{p+1}\omega_t^n dt  + \int_{t_0}^T \int_X \frac{\hat u^{p+1}}{p+1}( - R + \lambda + \tr_\omega\alpha   ) \omega_t^n dt.
\eea
Because of the lower bound of $R$ in Lemma \ref{lemma 1}, the last integral on the RHS is bounded above by $$C \int_{t_0}^T\int_X \frac{1}{p+1} \hat u^{p+1} \omega_t^n dt.$$  To deal with the term involving $\nabla_j \alpha_{\bar k i} \overline{S^k_{ij}}$ we make use of integration by parts and Cauchy-Schwarz inequality as follows
\bea\nonumber
\int_{t_0}^T \int_X 2 \xi \hat u^p Re(\nabla_j \alpha_{\bar k i}\overline{S^k_{ij}}) \omega_t^n dt & = & \int_{t_0}^T \int_X  - 2\xi \hat u^p Re( \alpha_{\bar k i} \overline{ S^k_{ij,\bar j}  }  ) - 2 \xi p \hat u^{p-1}Re ( \alpha_{\bar k i } \nabla_j \hat u  \overline{S^k_{ij}}   )\\
& \le & \int_{t_0}^T \int_X  C \xi \hat u^p |\bar \nabla S| +   Cp \xi \hat u^p ( |\nabla S| + |\bar \nabla S|  )\nonumber\\
&\le & \int_{t_0}^T \int_X C p^2 \xi \hat u^p + \frac{1}{2}\xi \hat u^p (|\nabla S|^2 + |\bar \nabla S|^2  ),\nonumber
\eea
where in the first inequality we apply Kato's inequality $|\nabla |S||\le |\nabla S| + |\bar \nabla S|$. Combining the estimates and varying $s\in [t_1, T']$, we conclude that
$$\sup_{t\in [t_1, T']} \int_X \hat u^{p+1}  \omega_t^n + \int_{t_1}^{T'} \int_X |\nabla \hat u^{\frac{p+1}{2}}  | ^2\omega_t ^n dt \le C \Big(  (p+1)^3 + \frac{1}{t_1 - t_0}   \Big) \int_{t_0}^{T'} \int_X \hat u^{p+1} \omega_t^n dt =:A. $$
We can now use the standard parabolic Moser iteration to conclude a sub-mean-value inequality(noting that the Sobolev inequality holds for $\omega_t$ by the assumption on the equivalence of the metrics  $\omega_t$ and $\omega_0$).

By H\"older inequality and the Sobolev inequality, we have
\bea\label{eqn:new 3}
\int_{t_1}^{T'} \int_X \hat u^{(p+1) (1+ \frac 1 n)} & \le & \int_{t_1} ^{T'} \Big( \int_X \hat u^{p+1}  \Big)^{1/n} \Big(  \int_X (\hat u^{\frac{p+1}{2}})^{ \frac{2n}{n-1}  }  \Big)^{(n-1)/n}dt\\
&\le  & \nonumber  A^{1/n} \int_{t_1}^{T'} C_0 \Big( \int_X  |\nabla \hat u^{\frac{p+1}{2}}|^2 + \hat u^{p+1}   \Big)\\
& \le & \nonumber 2 C_0 A^{(n+1)/n}.
\eea
Therefore
\bea\label{eqn:middle}
\Big( \int_{t_1}^{T'} \int_X \hat u^{\frac {n+1}{n} (p+1)} \Big) ^{\frac{n}{(n+1)(p+1)}} \le C_3^{1/(p+1)} \Big( (p+1)^3 + \frac{1}{t_1 - t_0}   \Big)^{1/(p+1)}
\Big(\int_{t_0}^{T'} \int_X \hat u^{p+1}  \Big)^{1/(p+1)} .\eea
If we denote $G(p,t) = \Big( \int_t ^{T'} \int_X \hat u^{p} \Big)^{1/p}$, the inequality (\ref{eqn:middle}) yields that
\bea\label{eqn:mid 1}
G(\frac{n+1}{n}(p+1) , t_1) \le C_3^{1/(p+1)} \Big( (p+1)^3 + \frac{1}{t_1 - t_0}   \Big)^{1/(p+1)} G(p+1, t_0).
\eea
For simplicity we denote $\eta = (n+1)/n>1$, and define a sequence of numbers $p_k + 1 = (p_{k-1} + 1) \eta =\cdots = (p_0 + 1) \eta^k$. For any $0< s_1 < s_2 < T'$, we define a sequence of times $t_k = s_1 + (1-\eta^{-k}) (s_2 - s_1)   $, $t_0 = s_1$ and $t_\infty = s_2$. Applying (\ref{eqn:mid 1}) for the pairs $(p_k+1, t_k)$, we get
$$ G ( p_{k+1} + 1, t_{k+1}  ) \le C_3^{\frac{1}{(p_0 + 1) \eta^k}} \Big(  (p_0+1)^3  + \frac{1}{s_2- s_1}  \Big)^{\frac{1}{(p_0+1)\eta^k}} \eta^{ \frac{3k}{(1+p_0) \eta^k}   } G(p_k + 1, t_k) ,  $$
iterating this inequality we get
\bea \nonumber G(\infty, s_2) & \le &C_3^{ \sum_k \frac{1}{(1+p_0) \eta^k} }  \Big(  (p_0+1)^3  + \frac{1}{s_2- s_1}  \Big)^{\sum_k \frac{1}{(p_0+1)\eta^k}} \eta^{ \sum_k \frac{3k}{(1+p_0) \eta^k}   } G(p_0 + 1 , s_1)\\
& \le &\nonumber  C_4 C_3^{ \frac{n+1}{1+p_0}}  \Big(  (p_0+1)^3  + \frac{1}{s_2- s_1}  \Big)^{ \frac{n+1}{(p_0+1)}} G(p_0 + 1 , s_1).
\eea
Let $p_0 = 1$, then the estimate above implies that
\bea\label{eqn:mid 2}
\sup_{X \times [s_2, T]} \hat u   \le C_5 ( 1+ \frac{1}{s_2 - s_1}   )^{(n+1)/2}\Big(\int_{s_1}^{T'} \int_X \hat u^2  \Big)^{1/2}.
\eea
Since we do not know the $L^4$-norm bound of $|S|$, we will use another  iteration argument to replace the $L^2$-norm of $\hat u$ on the RHS to the $L^1$-norm. We fix a $p \in (0,2)$, and denote $$h(s) = \sup_{X\times [s,T']} \hat u,$$ then (\ref{eqn:mid 2}) implies that
\bea\nonumber
h(s_2)& \le& C_5 (1+ \frac{1}{s_2 - s_1})^{(n+1)/2} h(s_1)^{\frac{2- p}{2}} \Big( \int_{s_1} ^{T'} \int_X \hat u^p \Big)^{1/2}  \\
&\le & \label{eqn:mid 4} \frac 1 2 h(s_1) + C_5^{p/2} (1+ \frac{1}{s_2 - s_1})^{(n+1)/p} \Big( \int_{s_1} ^{T'} \int_X \hat u^p \Big)^{1/p}.
\eea
Now for any $0\le t_1< t_2 \le T'$, we define a sequence of times $r_k = t_2 - (1-\delta^k) (t_2 - t_1)$, for some $\delta\in (0,1)$ to be determined later. Clearly $r_0 = t_2$ and $r_\infty = t_1$, and $r_k - r_{k+1} = (t_2 - t_1) (1-\delta) \delta^k$.     Iterating (\ref{eqn:mid 4}), we get
\bea\nonumber
h(t_2) = h(r_0) &\le& \frac 1 2 h(r_1) + C_5^{1/p} ( 1+ \frac{1}{r_0 - r_1}   )^{(n+1)/p} \Big(\int_{t_1}^{T'} \int_X \hat u^p  \Big)^{1/p}\\
& \le & \nonumber \frac{1}{2^k } h(r_k) +  C_5^{1/p} \Big(\int_{t_1}^{T'} \int_X \hat u^p  \Big)^{1/p} \sum_{i=0}^{k-1} 2^{-i} (1+ \frac{1}{r_i - r_{i+1}}   )^{(n+1)/p}\\
& \le & \nonumber \frac{1}{2^k } h(r_k) +  C_6^{1/p} \Big(\int_{t_1}^{T'} \int_X \hat u^p  \Big)^{1/p} \frac{\sum_{i=0}^{k-1} 2^{-i}  \delta^{-(n+1)i/p}}{(t_2 - t_1)^{(n+1)/p} (1-\delta)^{(n+1)/p}}.
\eea
If we choose $\delta\in (0,1)$ such that $2 \delta^{(n+1)/p}> 1$, then  the summation on the RHS of the above converges. Noting that $h$ is apriorily bounded (may not be uniform), and letting $k\to\infty$, we arrive at the desired estimate
\bea\nonumber
h(t_2)\le C_7^{1/p} \frac{1}{(t_2- t_1)^{(n+1)/p}} \Big( \int_{t_1}^{T'} \int_X \hat u^p   \Big)^{1/p}.
\eea
Setting $p = 1$, we get
$$\sup_{X\times [t_2 ,T']} \hat u\le \frac{C}{( t_2 - t_1)^{n+1} } \int_{t_1}^{t_2} \int_X \hat u \omega_t^n dt\le \frac{C}{(t_2- t_1)^{n+1}} \int_{t_1}^{T'} \int_X (|S|^2 + 1  ) \omega_t^n dt,$$
for any times $0\le t_1 < t_2 < T'<T$. It suffices to bound the $L^2$-norm of $|S|$ from above. By the Chern-Lu inequality (\cite{Y}), we have
$$(\frac{\partial }{\partial t} -  \Delta) \tr_{\omega_0} \omega \le \lambda \tr_{\omega_0} \omega + \tr_{\omega_0} \alpha  +  B \tr_\omega \omega_0 \tr_{\omega_0}\omega -  \delta_K |S|^2\le -  \delta_K |S|^2 + C_{K}  $$
where $\delta_K>0$ is a constant depending only on $K$ and the dimension $n$, and $-B$ is the lower bound of the bisectional curvature of $\omega_0$. Integrating the equation above  we obtain
\bea\nonumber
\frac{d}{dt}\int_X (\tr_{\omega_0}\omega) \omega^n &\le &  - \delta_K \int_X |S|^2 \omega^n + \int_X \tr_{\omega_0}\omega ( - R + \lambda + \tr_\omega\alpha  )\omega^n +C_{K}\\
&\le &\nonumber  - \delta_K \int_X |S|^2 \omega^n + C_{K}
\eea
where in the last inequality we use the lower bound of $R$ as in Lemma \ref{lemma 1}. Integrating over $t\in [ t_1, T']$ we get the desired $L^2$-bound of $|S|$ on $X\times [t_1 ,T']$. The proof of the upper bound $|S|^2$ is complete.

\medskip
Now that we have the third order estimates for the ptentials, we can come to the proof of the estimate for $\p_t\o$, which is a key estimate distinguishing the parabolic case from the elliptic case. Even so, we do not have as yet an estimate for the Riemann curvature tensor, which we circumvent below by a local Moser iteration argument:

\begin{lemma}\label{lemma new}
There is a constant $C=C(\omega_0,\kappa,K,\lambda)>0$ such that
$$\sup_X \left|\frac{\partial \omega_t}{\partial t}\right|\le C,\quad \forall ~t\in [0,T).$$
\end{lemma}
\noindent{\it Proof.}
Taking $\frac{\partial }{\partial t}$ on both sides of  the equation (\ref{eqn:varphi}) and adding the resulting equation to $\frac{1}{\kappa-1}\times$(\ref{eqn:F}), we get
$$\left(\frac{\partial}{\partial t} - \Delta\right) \left( \dot\varphi + \frac{\kappa}{\kappa-1} F \right) = \lambda\dot\varphi - \frac{\kappa}{\kappa-1}\tr_{\omega_t}\alpha_0 + \frac{\kappa}{\kappa-1} b.$$
We denote $\Phi = \dot\varphi + \frac{\kappa}{\kappa-1}F$, then we calculate (the norm of $\ddbar \Phi$ is under the metric $\omega = \omega_t$)
\bea\label{eqn:new 1}
&&(\frac{\partial}{\partial t} - \Delta)|\ddbar \Phi|^2\nonumber\\ & = & 2 \Phi_{\bar i j} \left(\lambda \Phi - \frac{\lambda\kappa}{\kappa-1} F - \frac{\kappa}{\kappa-1}\tr_{\omega_t} \alpha_0\right)_{\bar k l} g^{j\bar k} g^{l\bar i} - 2\lambda |\ddbar \Phi|^2 \\
&& \nonumber - 2\alpha_{\bar i j} \Phi_{\bar k l} \Phi_{\bar p q} g^{q \bar k } g^{j \bar p} g^{l\bar i} + \Phi^{j\bar i} \Phi^{l\bar k} R_{\bar i j \bar k l} - 2|\nabla \ddbar \Phi|^2.
\eea
Since the Riemannian curvature $Rm(\omega_t)$ on the RHS of (\ref{eqn:new 1}) is not a priori bounded, we cannot apply global (in space $X$) Moser iteration as in the proof of Lemma \ref{lemma 2} to bound $|\ddbar \Phi|^2$. Instead we will use the local expression of $Rm$, the Calabi $C^3$-estimate in Lemma \ref{lemma 2} and a local Moser iteration argument. To begin with, we can choose a cover of $X$ by Euclidean balls $\{B_a, z_a^i\}_a$, where $z^i_a$ are the complex coordinates. Without loss of generality we may assume $B_a$ are Euclidean balls with radius $3$ and the balls $\{\frac 1 3 B_a\}$ also cover $X$, here $\frac 1 3 B_a$ denotes the concentric ball with radius $1/3$ of that of $B_a$. Recall that the curvature $Rm$ is given by (in $B=B_a$)
$$R_{\bar i j\bar k l} =  - \frac{\partial ^2 g_{\bar i j}}{\partial z^l \partial z^{\bar k}} + g^{p\bar q} \frac{\partial g_{\bar q j}}{\partial z^l}\frac{\partial g_{\bar i p}}{\partial z^{\bar k}},$$
and the second term is uniformly bounded in $B$ by Lemma \ref{lemma 2}. We will use integration by parts to deal with the first term term in $R_{\bar i j \bar k l}$. We first observe that
$$\nabla_p \Phi_{\bar i j} = \frac{\partial ^3 \Phi}{\partial z^p \partial z^{\bar i}\partial z^j} - \Gamma^k_{pj} \Phi_{\bar i k},$$ $\Gamma^k_{pj}$ is bounded by Lemma \ref{lemma 2} so
$$| \nabla \ddbar \Phi  |^2 \ge \frac 1 2|D \ddbar \Phi|^2 - C |\ddbar \Phi|^2, $$ where we use $D$ to denote the ordinary derivatives in $(B, z^i)$.

We fix times $t_1< t_2 < T'< T$ and radii $r_2 < r_1 \le 3$, and define cut-off functions $\xi(t)$ as in the proof of Lemma \ref{lemma 2} and $\rho(z) = \rho(|z|)$ such that $\rho = 1$ on $B(0,r_2)$ and $\rho = 0$ outside $B(0,r_1)$, and $|D\rho|\le \frac{2}{r_1 - r_2}$. We denote $w = \max( |\ddbar \Phi|^2,1 )$. Multiply both sides of the equation (\ref{eqn:new 1}) by $\xi \rho^2 w^p$ (for $p\ge 1$) and do integration by parts, then we get
\bea\label{eqn:new 2}
&&\frac{d}{dt}\Big( \int_{B} \frac{\xi \rho^2 w^{p+1}}{p+1} \omega_t^n\Big) + \int_B \frac{p\xi}{(p+1)^2} |\nabla ( \rho w^{\frac{p+1}{2}}  )   |^2 - 2 \int_B \xi w^{p+1} |\nabla \rho|^2 \omega_t^n \\
& \le\nonumber & \int_B \frac{\xi \rho^2}{p+1} w^{p+1} (-R + \lambda n + \tr_{\omega_t} \alpha_t ) \omega_t^n + \int_B \frac{\rho^2 w^{p+1}\xi'(t)}{p+1} + \int_B C w^{p+1}\\
&& \nonumber - \frac 1 2 |D\ddbar \Phi|^2 -\frac{2\kappa\xi}{\kappa-1}\rho^2 w^p \Phi^{j\bar i } (\lambda F + \tr_{\omega_t} \alpha_0)_{\bar i j} - \xi \rho^2 w^p \Phi^{j\bar i} \Phi^{k\bar l} \frac{\partial}{\partial z^k} \Big( \frac{\partial g_{\bar ij }}{\partial z^{\bar l}} \Big),
\eea
we deal with the last two terms in the integral by integration by parts as follows: the second last term is equal to
\bea\nonumber
&& \int_B \frac{2\kappa\xi}{\kappa-1}\rho^2 ( 2 \rho^{-1} w^p \rho_j \Phi^{j\bar i  }  + p w^{p-1}  w_j   \Phi^{j\bar i} +   w^p \partial_j \Phi^{j\bar i} +  w^p \Phi^{j\bar i} \partial_j \log \det g )(\lambda F + \tr_{\omega_t} \alpha_0)_{\bar i} \\
&\le & \int_B C p^2 \xi \rho^2 w^{p+1} + C \xi w^{p+1} |\nabla \rho|^2 + \frac{1}{10} \xi \rho^2 w^p |D \ddbar\Phi|^2,\nonumber
\eea
where in the inequality above we use the known estimates that $|\nabla F|\le C$, $|\nabla \tr_{\omega_t} \alpha_0|\le C$ and $|\nabla \dot\varphi|\le C$ which follows from the equation (\ref{eqn:varphi}) and Calabi $C^3$-estimates. By integration by parts the last integral in (\ref{eqn:new 2})  is equal to
\bea
\nonumber && \int_B \xi \rho^2( \rho^{-1} w^p \rho_k \Phi^{j\bar i} \Phi^{k\bar l} +  p w^{p-1} w_k \Phi^{j\bar i}\Phi^{k\bar l} + w^p \partial_k (\Phi^{j\bar i}\Phi^{k\bar l}   ) + w^p  \Phi^{j\bar i}\Phi^{k\bar l} \partial_k \log \det g     ) \frac{\partial g_{\bar i j}}{\partial z^{\bar l}}\\
&\le & \int_B C p^2 \xi \rho^2 w^{p+1} + C \xi w^{p+1} |\nabla \rho|^2 + \frac{1}{10} \xi \rho^2 w^p |D \ddbar\Phi|^2.\nonumber
\eea
Substituting the above inequalities to (\ref{eqn:new 2}) and argue as in the proof of Lemma \ref{lemma 2}, we get
\bea\nonumber
&& \sup_{[t_2, T']}\int_{B_{r_2}} w^{p+1}\omega_t^n + \int_{t_2}^{T'} \int_{B_{r_2}} |\nabla (\rho w^{\frac{p+1}{2}})  |^2 \omega_t dt\\
&\le & C \Big( (p+1)^3 + \frac{1}{(r_1 - r_2)^2} + \frac{1}{- t_1 +t_2}\Big) \int_{t_1}^{T'} \int_{B_{r_1}} w^{p+1}\omega_t^n dt.
\eea
By Holder inequality and Sobolev inequality as in (\ref{eqn:new 3}), we have $$ H ( (p+1)\eta, t_2, r_2 )\le C^{1/(p+1)} \Big( (p+1)^3 + \frac{1}{(r_1 - r_2)^2} + \frac{1}{- t_1 +t_2}\Big)^{1/(p+1)} H(p+1,t_1,r_1)   $$
where $H(p, t', r) =\Big(\int_{t'}^{T'} \int_{B_r} w^p \omega_t ^n dt  \Big)^{1/p} $. Taking $p_{k} + 1 = \eta^k (p_0+1)$, $r_k = 1+ \eta^{-k/2}$ and $t_k = T' - 1 - \eta^{-k}$, we can do the iteration argument and the trick as in the proof of Lemma \ref{lemma 2} to conclude that
$$\sup_{[T'-1,T']\times B_1} w\le \int_{T'-2}^{T'} \int_{B_2} w \, \omega_t^n dt.    $$
As the last step, we have the following equation from Bochner formula
\bea\nonumber(\frac{\partial}{\partial t} - \Delta)|\nabla\Phi|^2 & = &- |\ddbar \Phi|^2 - |\nabla \nabla \Phi|^2 - \lambda |\nabla \Phi|^2 - \alpha(\nabla \Phi,\bar \nabla \Phi) \\
&& \nonumber- 2 Re \< \nabla \Phi, \bar \nabla ( -\dot\varphi - \frac{\kappa}{\kappa - 1}\tr_{\omega}\alpha_0   ) \> \le  - |\ddbar \Phi|^2 + C,
\eea
integrating this equation over $[T'-2,T']\times X$, and by the lower bound of $-R$ we can obtain the $L^2( [T'-2,T']\times X  )$-norm bound of $|\ddbar \Phi|$, hence that of $w$. This completes the proof of the bound of $|\ddbar \Phi|$ since the balls $B_1$ cover $X$ by the original choice. Since $\ddbar F$ is uniformly bounded by assumption, this implies that $\ddbar \dot\varphi$ is bounded, so is $\frac{\partial \omega_t}{\partial t}$. The proof of Lemma \ref{lemma new} is complete.

\bigskip

We apply the linear theory of parabolic equations, which requires the coefficients $a^{ij}(x,t)$ in (\ref{eqn:linear equation}) below to be continuous with uniform modulus of continuity or uniform H\"older regularity in {\bf both space and time directions}. Lemma \ref{lemma 2} gives the uniform Lipschitz continuity of $g^{i\bar j}$ in space directions and Lemma \ref{lemma new} provides the Lipschitz regularity of $g^{i\bar j}$ in the time direction.

\bigskip

To continue with the higher order estimates, we use the standard linear theory for parabolic equations. For convenience we state the following $W^{2,1}_{p}$-estimate (see Theorem 7.22 in \cite{Lie}). For a smooth bounded domain $\Omega\subset{\mathbf R}^m$ and $T>0$ we denote the parabolic cylinder $$\Omega_T:= \Omega\times [0, T]\subset{\mathbf R}^{m+1},$$ and  the $W^{2,1}_p(\Omega_T)$-norm of a function $u$ in $\Omega_T$ is defined as
$$\| u\|_{W^{2,1}_{p}(\Omega_T)} = \sum_{j=0}^2 \| D^j u\|_{L^p(\Omega_T)} + \| \partial_t u\|_{L^p(\Omega_T)}.   $$
\begin{lemma}\label{Lp estimate}
Let $u\in W^{2,1}_{p, loc}(\Omega_T)\cap L^p(\Omega_T)$ for some $p\ge 2$ satisfy the parabolic equation
\bea\label{eqn:linear equation}\frac{\partial u}{\partial t}(x,t) = a^{ij}(x,t) u_{ij}(x,t) + f(x,t), \quad{\mathrm{in}}~ \Omega_T\eea where $a^{ij}$ is strictly elliptic in the sense that $\lambda |\xi|^2 \le a^{ij} \xi_i \xi _j \le \Lambda |\xi|^2$ for some $\lambda, \Lambda >0$ and any $\xi\in{\mathbf R}^m$ and $a^{ij}$ is continuous in $\Omega_T$. Suppose $f\in L^p(\Omega_T)$, then for any domain $\Omega'\subset\subset \Omega$ and $\varepsilon_0>0$ there exists a constant $C = C(\lambda ,\Lambda, d(\Omega', \partial\Omega), \omega_a,  \varepsilon_0, m, p)>0$ such that
$$ \| D^2 u\|_{L^p( \Omega'\times [\varepsilon_0,T]  )} + \| \frac{\partial u}{\partial t} \|_{L^p( \Omega'\times [\varepsilon_0,T]  ))} \le C ( \| f\|_{L^p(\Omega_T)} + \| u\|_{L^p(\Omega_T)}      ),  $$ where $\omega_a$ denotes the modulus of continuity of $a^{ij}$.
\end{lemma}
In particular, if $f\in L^\infty(\Omega_T)$ and $u\in L^\infty(\Omega_T)$,  by Sobolev embedding, we have $u\in C^{ 1+\beta, \frac{1+\beta}{2}  }( \Omega_T\backslash \partial_P \Omega_T )$ for any $\beta\in (0,1)$. Here $\partial_P \Omega_T$ denotes the parabolic boundary of $\Omega_T$. It follows that we have the uniform estimate
\bea\nonumber
\| u\|_{C^{1+\beta, \frac{1+\beta}{2}}( \Omega'\times[\varepsilon_0, T]  )} \le C ( \| f\|_{L^\infty(\Omega_T)} + \| u\|_{L^\infty(\Omega_T)}  ),
\eea
where the constant $C$ depends in addition on $\beta$.

\medskip

We make the following elementary observation from elliptic theory.
\begin{lemma}\label{lemma 3}
Under the assumptions (\ref{eqn:assumption}), there is a constant $C=C(K,\omega_0)>0$ such that
$$\sup_{X\times [0,T)}(|\nabla F| + |\nabla \varphi|) \le C. $$
\end{lemma}
\noindent{\it Proof.}
By the assumption (\ref{eqn:assumption}), we have
$$| \Delta_{\omega_0} ( \varphi - \bar \varphi  )  | \le C, \quad |\Delta_{\omega_0} ( F - \bar F)| \le C,$$
where $\bar \varphi$ and $\bar F$ denote the average of $\varphi$ and $F$ with respect to $\omega_0^n$, respectively. It is elementary to see the $C^0$ bound of $\varphi - \bar \varphi$ and $F-\bar F$. The desired gradient bound follows standard linear theory for elliptic equations. In fact we can also get the $C^\beta$-bound of $\nabla\varphi$ and $\nabla F$ for any $\beta\in (0,1)$, though this is not needed for later calculations.
\hfill \qed

\medskip

We choose (and fix) finite coordinate charts of $X$, $\{U_a\}$, $\{ U_a'\}$ and $\{U_a''\}$ such that $U_a''\subset\subset U_a'\subset \subset U_a$ and $X = \cup_a U_a''$, and $\{z^i_{a}\}$ the complex coordinate functions on $U_a$. In the following we will work on any chosen $U_a$ and omit the subscript $a$ in $U_a$, $U_a'$, $U_a''$ and $z^i_a$. We will fix times $t_0\ge 0$, $t_1 = t_0+1$ and $t_2 = t_0 + 2$.

\medskip

\begin{lemma}\label{lemma 5}
For any $\beta\in (0,1)$, there exists a constant $C=C(\omega_0, \kappa, K, \beta)>0$ such that $$\| \ddb \varphi\|_{C^{1+\beta, \frac{1+\beta}{2}}(X\times [t_1,t_2])}\le C$$ or equivalently
$$\| \omega_t\|_{C^{1+\beta, \frac{1+\beta}{2}} (X\times [t_1,t_2])  }\le C.$$
\end{lemma}
\noindent{\it Proof.}
We work on a coordinate chart $U=U_a$ as before. Taking $\frac{\partial^2}{\partial z^i\bar \partial z^j}$ on both sides of (\ref{eqn:varphi}) we get on $U\times [t_0,t_2]$
\bea\nonumber \frac{\partial \varphi_{\bar ji}}{\partial t}&  =  & \Delta_{\omega_t} \varphi_{\bar ji} - g^{p\bar m }\partial_{\bar j} g_{\bar m k} g^{k \bar q} \varphi_{i\bar q p} + g^{p\bar q} \partial^2_{\bar j i} \hat g_{\bar q p} - g^{p\bar m} \partial_{\bar j} g_{\bar m k } g^{k\bar q} \partial_i \hat g_{\bar q p} + \hat R_{\bar j i}\\
&& - h_{0, \bar j i} - \varphi_{\bar j i} - F_{\bar j i} =: \Delta_{\omega_t} \varphi_{\bar j i} + f_1\label{eqn:varphi ij}
\eea
where $f_1\in L^\infty( U\times[t_0,t_2]  )$ by assumptions and Lemma \ref{lemma 2}. Applying the linear theory  in Lemma  \ref{Lp estimate} and Sobolev embedding theorem,  we conclude that for any $\beta\in (0,1)$ there is a uniform $C=C(\kappa,K ,\beta,\omega_0)>0$ such that  ($\hat t_1 = t_1-1/2$ and $p = \frac{2n+2}{1-\beta}$)
\bea\nonumber
 \| \varphi_{\bar j i} \|_{C^{1+\beta,\frac{1+\beta}{2}}( U'\times [t_1,t_2]  )} &\le & \| D^2\varphi_{\bar j i} \|_{L^p( U'\times [\hat t_1,t_2]  )} + \| \partial_t\varphi_{\bar j i} \|_{L^p( U'\times [\hat t_1,t_2]  )} \\
& \le &  C(\| f_1\|_{L^\infty( U\times [t_0,t_2]  )}  + \| \varphi_{\bar j i}\|_{L^\infty(U\times [t_0,t_2])} )\le C. \label{eqn:new 6}
\eea
Since $\{U_a'\}$ covers $X$, the lemma follows.
\hfill \qed

\begin{lemma}\label{lemma 4}
For any $\beta\in (0,1)$, there exists a constant $C=C(\omega_0, \kappa, K, \beta)>0$ such that $$\| \nabla F\|_{C^{1+\beta, \frac{1+\beta}{2}} (X\times [t_1,t_2])   }\le C, \quad{\mathrm{and}~  } ~\| \alpha\|_{C^{\beta,\frac{\beta}{2}}(X\times [t_1,t_2])}\le C   .$$
\end{lemma}
\noindent{\it Proof.}
We work on a coordinate chart $U = U_a$. Taking $\frac{\partial }{\partial z^i}$ on both sides of the equation (\ref{eqn:F}) we have the following equation holds on $U\times [t_0,t_2]$
\bea\label{eqn:Fi} \frac{\partial F_i}{\partial t} = \kappa \Delta_{\omega_t} F_i - \kappa g^{p\bar m} \partial_i g_{\bar m k} g^{k\bar q} F_{\bar q p} - \kappa \partial_i (\tr_{\omega_t}\alpha_0)=: \kappa \Delta_{\omega_t} F_i + f.    \eea
By the Calabi $C^3$-estimates in Lemma \ref{lemma 2} we see that $\| f\|_{L^\infty(U\times [t_0,t_2]  )}\le C$ for a uniform constant $C>0$.
Since the coefficients of $\kappa\Delta_{\omega_t}$ is Lipschtiz continuous with a uniform Lipschtiz constant by Lemmas \ref{lemma 2} and \ref{lemma new},   we can now apply the linear theory in Lemma \ref{Lp estimate} to conclude that for any $p\ge 2$
\bea\label{eqn:new 5}
&& \| D^2 F_i\|_{L^p(U'\times[t_1-1/2, t_2])} + \| \partial_t F_i\|_{L^p(U'\times[t_1-1/2, t_2])} \\ & \le& C(p) (   \| f\|_{L^\infty(U\times [t_0,t_2])} + \| F_i\|_{L^\infty(U\times [t_0,t_2])}   )\le C,\nonumber
\eea thanks to Lemma \ref{lemma 3} on the gradient bound of $F$.
By Sobolev embedding then we get
$\| F_i \|_{C^{1+\beta, \frac{1+\beta}{2}} (U'\times [t_1,t_2])} \le C.$ Taking $\frac{\partial }{\partial z^{\bar j}}$ on both sides of (\ref{eqn:Fi}) we get the equation for $F_{\bar j i}$
\bea\label{eqn:Fij}
\frac{\partial F_{\bar j i}}{\partial t} = \kappa \Delta_{\omega_t} F_{\bar j i } + \kappa \partial_{{\bar j}}g^{m\bar k} F_{i \bar k  m} + f_{\bar j}= : \kappa \Delta_{\omega_t} F_{\bar j i} + f_2.
\eea
We observe that the terms in $f_{2}$ in (\ref{eqn:Fij}) involve either $\partial^{2}_{\bar k l }g_{\bar p q}$ or $D F_{\bar k l}$, besides the bounded factors. By (\ref{eqn:new 6}) and (\ref{eqn:new 5}) we know $\partial^{2}_{\bar k l }g_{\bar p q}$ and $D F_{\bar k l}$ are both bounded in $L^p(U'\times[t_1-1/2, t_2]  )$ for any $p\ge 2$, so is $f_2$ by Holder inequality. Applying linear theory in Lemma \ref{Lp estimate} and Sobolev embedding theorem again, we conclude that ($\hat t_1 = t_1-1/2$ and $p = 2\cdot \frac{2n+2}{1-\beta}$)
\bea\nonumber
\|F_{\bar j i}\|_{C^{1+\beta,\frac{1+\beta}{2}} ( U''\times[t_1,t_2]  )  } & \le & C ( \| D^2 F_{\bar j i}\|_{L^p(U'\times[ t_1,t_2])} + \| \partial_t F_{\bar j i}\|_{L^p(U'\times [t_1,t_2])}    )\\
&\le \nonumber& C ( \| F_{\bar j i}\|_{L^p(U\times [\hat t_1,t_2])} + \| f_2\|_{L^p(U\times [\hat t_1,t_2])}     )\le C.
\eea
 Since $\{ U_a''\}$ covers $X$, the lemma follows from the fact that $\alpha_t = \alpha_0 + i\ddb F$.

\bigskip


\noindent{\it Proof of Theorem \ref{higherorder}.} We can now apply the standard bootstrap argument to finish the proof of Theorem \ref{higherorder}. Taking $\frac{\partial}{\partial z^{\bar k}}$ on both sides of (\ref{eqn:varphi ij}), the lower order terms in the resulted equation are bounded uniformly in $L^p$ for any $p\ge 2$ by Lemmas \ref{lemma 5} and \ref{lemma 4}, and this gives rise to $W^{2,1}_{p, loc}$ estimates of $\varphi_{\bar j i \bar k}$. Feed these estimates into the equation obtained by taking $\frac{\partial}{\partial z_{\bar k}}$ on both sides of (\ref{eqn:Fij}). By similar argument this yields $W^{2,1}_{p, loc}$ estimates of $F_{\bar j i \bar k}$. We can repeat this process any finite number of times to conclude that $D^k \varphi_{\bar j i}$ and $D^k F_{\bar j i}$ are both in $W^{2,1}_{p, loc}$ for any $p\ge 2$ and $k\in{\mathbf Z}_+$. The higher order estimates of $\varphi_{\bar j i}$ and $F_{\bar j i}$ then follows from Sobolev embedding theorem, so do the higher order derivative estimates of $\omega_t$ and $\alpha_t$. Theorem \ref{higherorder} is proved.

\bigskip
We remark that perhaps the most important consequence of the assumption (\ref{eqn:assumption}) is that the Sobolev constant is uniformly bounded along the flow. This is known to be true, but highly non-trivial, along the K\"ahler-Ricci flow \cite{P, Zh}. It would clearly be desirable to determine whether, or under what conditions, the Sobolev constant will be uniform along the $\kappa$-LYZ flow.


\section{The Case of Riemann Surfaces}
\setcounter{equation}{0}

Both the Ricci flow and the Calabi flow have been completely solved on Riemann surfaces \cite{H, Ch, Str}. Thus the case of Riemann surfaces provides an excellent laboratory where to examine the complications arising from couplings to an additional field. It is easy to see that stationary points for the $\kappa$-LYZ flow always exist. Even so, the convergence of the $\kappa$-LYZ flow is not evident. What we can prove is the following:

\begin{theorem}
\label{RS}
Consider the $\kappa$-LYZ flow on a compact Riemann surface $X$ with negative Euler characteristic. Write $\alpha_t=\tau_t\o_t$ for a smooth function $\tau_t$. Assume that at initial time, we have
\bea
\label{nec3}
R(\omega_0)<0 \ {\rm and}\ \tau_0>0
\eea
everywhere on $X$. Then there exists a positive constant $\kappa_0$ such that for any $\kappa>\kappa_0$, $\kappa\not=0$, the $\kappa$-LYZ flow converges in $C^\infty$ to a pair $(\o_\infty,\alpha_\infty)$ with
\bea
R(\omega_\infty)=\lambda+\tau_\infty,
\qquad
\tau_\infty=\int_X\alpha_0
\eea
both constants. Actually, $\kappa_0$ can be taken to be any constant with $R_0+{1\over\kappa_0-1}\tau_0<0$.
\end{theorem}

\bigskip
We observe that the conditions on the initial data of the theorem imply that $\lambda<0$. Assuming this, it is easy to construct many initial data satisfying the conditions. Given a negatively curved compact Riemann surface $(X,\omega_0)$, there exists a sufficiently negative constant $\lambda$ such that the cohomology class $-\lambda[\omega_0]+2\pi c_1(X)$ lies in the K\"ahler cone. We can then choose a positive closed $(1,1)$-form $\alpha_0$ in this class.

\medskip
\noindent
{\it Proof of Theorem \ref{RS}}. The $\kappa$-LYZ flow had been written
earlier in terms of potentials (\ref{eqn:varphi}) and (\ref{eqn:F})). Here, we shall exploit the setting of Riemann surfaces in order to write the flow in terms of the conformal factor $\phi$ and the trace $\tau$ defined by
\bea
\omega=e^{\phi}\omega_0,
\qquad \tau\equiv {\rm tr}_\o\alpha=g^{z\bar z}\alpha_{\bar zz}.
\eea
Clearly the flow (\ref{LYZo}) can be rewritten in terms of $\phi$ as
\bea
\dot\phi = \frac{\partial\phi}{\partial t}= - R+\lambda+\tau,
\eea
while the flow (\ref{LYZa-theta}) for $\alpha$ can be rewritten in terms of $\tau$ as
\bea
\label{eqn:tau}
\p_t\tau= \kappa\Delta\tau-\tau(-R+\lambda+\tau).
\eea
We can rewrite the equation (\ref{eqn:for scalar}) for the evolution equation of the scalar curvature $R$ as follows
\bea
\label{eqn:r}
\p_tR
&=&\Delta R-\Delta\tau-R(-R+\lambda+\tau).
\eea
We  have seen the advantage of taking the coefficient $\kappa \neq 1$ in (\ref{LYZa-theta}) in obtaining the heat inequality (\ref{kappa}). For Riemann surfaces, this inequality is strengthened into a heat equation for $R+\frac{1}{\kappa-1}\tau$,
\bea
\label{eqn:r and tau}
\p_t\left(R+\frac{1}{\kappa-1}\tau\right)
=
\Delta\left(R+\frac{1}{\kappa-1}\tau\right)-\left(R+\frac{1}{\kappa-1}\tau\right)(-R+\lambda+\tau).
\eea

\subsection{Estimates for $\tau$ and $R$}

From the equations (\ref{eqn:tau}) and (\ref{eqn:r and tau}), it follows by the maximum principle that
$$ \max_X \left( R(\omega_t) + \frac{1}{\kappa-1}\tau_t  \right) <0, \quad {\mathrm{and}}\quad \min_X \tau_t >0,  $$ therefore $\max _X R(\omega_t) <0$. From (\ref{eqn:tau}) we have  at the maximum point of $\tau_t$
\bea
&&\frac{d \tau_{\max}}{dt}\le  \tau_{\max} R - \tau_{\max} \lambda - \tau^2_{\max}\nonumber\\ &\le& -\frac{\kappa}{\kappa-1}\tau^2 _{\max} - \lambda \tau_{\max}= -\frac{\kappa}{\kappa-1}\tau_{\max} \left( \tau_{\max} + \frac{\kappa-1}{\kappa}\lambda \right),\label{eqn:eqn for tau}
\eea
\noindent $\bullet$ if $\tau_{\max}|_{t= 0} < (\kappa-1)|\lambda|/\kappa$,  by (\ref{eqn:eqn for tau}), it follows that $\tau_{\max} + (\kappa-1)\lambda/\kappa < 0 $ for all $t>0$.

\medskip

\noindent $\bullet$ if $\tau_{\max}|_{t = 0} \ge (\kappa-1)|\lambda|/\kappa$, by solving the ODE
$$f' = -\frac{\kappa}{\kappa-1}f^2 - \lambda f,\quad f( 0 ) = \tau_{\max}|_{t=0} = : a\ge (\kappa-1)|\lambda|/\kappa,$$ and comparing $\tau_{\max}$ with this $f$ we get
$$\tau_{\max}(t) \le f(t) = (\kappa-1)\frac{|\lambda|}{\kappa} \frac{  \frac{a}{a+ (\kappa-1)\lambda/\kappa}   } {\frac{a}{a+(\kappa-1)\lambda/\kappa}  - e^{\lambda t}  } \le (\kappa-1)\frac{|\lambda|}{\kappa} \frac{  \frac{a}{a+ (\kappa-1)\lambda/\kappa}   } {\frac{a}{a+(\kappa-1)\lambda/\kappa}  -1  } .$$
Therefore we get in both cases $\tau_{\max} \le C(a, \lambda)$.

Similarly we deal with $R+\frac{1}{\kappa-1}\tau$. Applying maximum principle to $R+\frac{1}{\kappa-1}\tau$ in (\ref{eqn:r and tau}) we get
$$ \frac{d}{dt} \left(R+\frac{1}{\kappa-1}\tau\right)_{\min} \ge - \left(R+\frac{1}{\kappa-1}\tau\right) _{\min}\Big ( - \left(R+\frac{1}{\kappa-1}\tau\right)_{\min} + \lambda + \frac{\kappa}{\kappa-1} \tau     \Big).   $$
For notation simplicity we denote $h = - (R+\frac{1}{\kappa-1}\tau)_{\min}>0$,  then
$$\frac{d h}{dt} \le - h^2 - h \left(\lambda+\frac{\kappa}{\kappa-1}\tau\right) \le - h^2 + |\lambda| h = - h ( h -  |\lambda| ),$$
by similar argument as above, if $h|_{t= 0} < |\lambda|$, then $h\le |\lambda|$ for all $t>0$; otherwise we have
$$h(t) \le |\lambda| \frac{\frac{b'}{b' - |\lambda|}}{\frac{b'}{b'-|\lambda|} - 1} = C(\lambda, b'),\quad {\mathrm{if}}\quad h_{t= 0 }:= b' \ge |\lambda|.$$
So we find that $ - C(\lambda, b') \le R+\frac{1}{\kappa-1}\tau <0 $ for all $t>0$. On the other hand, from $0< \tau \le C(\lambda, a)$, we derive that  the scalar curvature $R$ satisfies $ - C(\lambda, a, b') \le R <0  $. Let us summarize what we have proven so far:

\begin{lemma}\label{prop:1}
Under the conditions in Theorem  \ref{RS} on the initial values, there is a uniform constant $C$ depending only on the initial values, $\lambda$ and $\kappa>\kappa_0$ such that
$$0< \tau \le C,\qquad  - C \le R <0.$$
\end{lemma}

\subsection{A Liouville-Entropy functional}

Motivated by the Liouville energy studied in conformal geometry (see for example \cite{Str}), we define a functional $E(\phi, \tau)$ associated to the LYZ flow as follows, assuming that $\tau$ is a positive function:
\bea
\label{energy}
E(\phi, \tau) = \int_X (\frac{1}{2} |\nabla \phi|_{\omega_0}^2 + R_0 \phi - \lambda e^{\phi} - \tau e^{\phi}) \omega_0 + \int_X (\tau \log \tau) e^{\phi}\omega_0,
\eea
where $R_0 = R(\omega_0)$ is the scalar curvature of $\omega_0$. Clearly $E$ is the sum of an ``energy part'' and an ``entropy part''. We have

\begin{lemma}
\label{energylemma}
We consider the $\kappa$-LYZ flow (\ref{LYZo},\ref{LYZa-theta}). Then

{\rm (a)} The functional $E(\phi,\tau)$ is monotone non-increasing along the flow.

{\rm (b)} Under the assumptions of Theorem \ref{RS} on the initial data $(\o_0,\alpha_0)$ and  the choice of $\kappa>\kappa_0$, there is a constant $C>0$ such that
\bea\label{eqn:lower energy}
E(\phi, \tau)\ge \int_X \frac{1}{4} |\nabla \phi|_{\omega_0}^2\;  \omega_0
+
2\pi|\chi|\cdot |\bar\phi|- C,\quad \forall ~ t\ge 0.
\eea
where $\chi$ is the Euler characteristic of $X$, and $\bar\phi$ is the average of the function $\phi$ with respect to the initial metric $\o_0$.

{\rm (c)} For any $p\geq 1$, there exists a constant $C_p$ so that, along the $\kappa$-LYZ flow, we have
\bea
\label{Lp}
\int_X e^{p|\phi|}\o_0\leq C_p.
\eea
\end{lemma}

\noindent
{\it Proof.} To prove (a), we
calculate the variation of $E$ along the flow: (we denote $\dot\phi = \frac{\partial \phi}{\partial t}$)
\bea
\frac{d}{dt} E(\phi,\tau) & = & \int_X  ( -\Delta_{\omega_0} \phi + R_0 - \lambda e^{\phi} - \tau e^\phi    ) \dot \phi \omega_0 - \int_X \dot \tau e^{\phi}\omega_0 \nonumber\\
&  & + \int_X( (\dot \tau \log \tau) e^\phi \omega_0 + \dot \tau e^\phi \omega_0 + (\tau\log \tau ) \dot\phi e^\phi\omega_0 \nonumber    )\\
& = & \int_X -(\dot \phi )^2 e^\phi \omega_0 + \int_X (  \kappa e^{-\phi} \Delta_{\omega_0} \tau  ( \log \tau  ) e^{\phi} \omega_0 \nonumber - \tau \dot \phi (\log \tau) e^{\phi} \omega_0 + (\tau \log\tau) \dot \phi e^{\phi}\omega_0      )\nonumber\\
& = & -\int_X (\dot\phi)^2 e^{\phi} \omega_0  - \kappa \int_X \frac{|\nabla \tau|_{\omega_0}^2}{\tau} \omega_0\nonumber \\
& = & -\int_X (\dot\phi)^2  \omega_t  - \kappa \int_X \frac{|\nabla \tau|_{\omega_t}^2}{\tau} \omega_t \le 0.\nonumber
\eea
This proves (a).

To prove (b), we note that the volume of $\o_t$ has been normalized to be $1$, and thus
\bea
\label{normalization}
\int_X e^\phi\o_0=1.
\eea
Since $\tau\log\tau$ is bounded from below, it follows that the contribution of the entropy term is bounded from below. Furthermore, we have seen that under the conditions of the lemma, $\lambda$ is negative and $\tau$ is bounded. Thus we can write
\bea
E(\phi,\tau)
&\geq&
\int_X{1\over 2}|\nabla_0\phi|^2\o_0
+R_0\phi\o_0-C\nonumber\\
&=&
\int_X{1\over 2}|\nabla_0\phi|^2\o_0
+R_0(\phi-\bar\phi)\o_0+2\pi\chi\bar\phi-C
\geq \int_X{1\over 4}|\nabla_0\phi|^2\o_0+2\pi \chi\bar\phi-C
\nonumber
\eea
where we changed constants and applied the Poincare inequality to absorb the integral of $R_0(\phi-\bar\phi)$ into the first term on the right hand side. Applying Jensen's inequality to the normalization equation (\ref{normalization}) shows that $\bar\phi\leq 0$. Since $\chi<0$, we obtain (b).

Finally, (c) follows from (b) and the Trudinger inequality in two dimensions, since (b) implies that both $\|\nabla\phi\|_{\o_0}$ and $\bar\phi$ are bounded. The proof of the lemma is complete.

\subsection{Uniform bounds for $\phi$}
\label{C00}

With the uniform bound for $R$ and the $L^p$ estimates for $e^{\phi}$ in (\ref{Lp}), we can apply the arguments of Struwe \cite{Str} to obtain uniform bounds
for $\phi$,
\bea
\label{C0}
\|\phi\|_{C^0}\leq C.
\eea
The details are as follows. First, we state the concentration/compactness theorem due to Struwe (Theorem 3.2, \cite{Str}) as a lemma.

\begin{lemma}\label{struwe}
Let $\o_j=e^{\phi_j}\o_0$ be a sequence of metrics on a compact Riemann surface $X$, with $[\o_j]=1$ and Calabi energy $\|R_j\|_{L^2(\o_j)}^2$ uniformly bounded. Then

{\rm (a)} Either there exists a constant $C$ so that $\|\phi_j\|_{H_{(2)}(\o_0)}\leq C$
for all $j$;

{\rm (b)} or there exists a finite set $Z$ of points $z_1,\cdots, z_k$
such that, for any $z_k\in Z$ and any neighborhood $V$ around $z_k$, we have
\bea
\label{concentration}
{\rm liminf}_{j\to\infty}\int_V|R_j|\o_j\geq 2\pi.
\eea
The cardinality of $Z$ can be bounded in terms of the uniform upper bound on the Calabi energy.
\end{lemma}

Returning to the proof of the uniform bound for $\phi$, it suffices now to show that the $L^p$ bound for $e^{|\phi|}$ and the uniform bound for $R$ rule out the existence of points $z_k$ satisfying the condition (\ref{concentration}). The above lemma of Struwe would imply then that
$\|u_j\|_{H_{(2)}(\o_j)}\leq C$, which impliies, in two real dimensions, that
there exists a positive constant $\alpha\in (0,1)$ with $\|\phi\|_{C^\alpha}\leq C$.

For each $z$ in $X$ and each $r>0$, let $B(z,r)$ be the ball centered at $z$ and of radius $r$ with respect to the fixed metric $\o_0$. Then we write
\bea
\int_{B(z,r)}|R_j|\o_j
\leq
\Big(\int_{B(z,r)}|R_j|^2\o_j\Big)^{1\over 2}\Big(\int_{B(z,r)}\o_j\Big)^{1\over 2}
\leq A \Big(\int_{B(z,r)}\o_j\Big)^{1\over 2}
\eea
where $A$ is the uniform upper bound for the Calabi energy. Next,
\bea
\int_{B(z,r)}\o_j
&=&
\int_{B(z,r)}e^{\phi_j}\o_0
\leq
\Big(\int_{B(z,r)}e^{2\phi_j}\o_0\Big)^{1\over 2}
\Big(\int_{B(z,r)}\o_0\Big)^{1\over 2}
\nonumber\\
&\leq &
\Big(\int_{B(z,r)}e^{2\phi_j}\o_0\Big)^{1\over 2}.
\eea
But the uniform bound for the integral of $e^{2|\phi_j|}$ obtained in (\ref{Lp}) shows that the right hand side tends to $0$ as $r\to 0$. This rules out the concentration of the curvature, and the proof of the $C^\alpha$ bounds for $\phi_j$ is complete. \hfill\qed

\subsection{Convergence of $\|\nabla \tau\|_{L^2}^2$ and $\|\dot\phi\|_{L^2}^2$ to $0$}

The bound from below for the functional $E(\phi,\tau)$ implies the existence of a sequence $t_j$, $j\leq t_j\leq j+1$ with
\bea
\int_X (\dot\phi)^2+\int_X |\na\tau|^2
\to 0
\eea
as $j\to \infty$. (Here the integrals are with respect to the metric $\o_{t_j}$, which we do not write explicitly for notational simplicity.) Indeed, for any time $T$, it follows from the formula for the time derivative of $E(\phi,\tau,t)$ that
\bea
E(T) -E(0)=-\int_0^T\int_X
(\dot\phi)^2  \omega_t  + \kappa \int_X \frac{|\nabla \tau|_{\omega_t}^2}{\tau} \omega_t
dt\eea
and hence, applying the lower bound for $E(T)$ and letting $T\to\infty$,
\bea
\int_0^\infty\bigg(\int_X
(\dot\phi)^2  \omega_t  + \kappa \int_X \frac{|\nabla \tau|_{\omega_t}^2}{\tau} \omega_t\bigg)
\,dt<\infty.
\eea
Since $1/\tau$ is uniformly bounded from below, the convergence of the integral in $t$ implies the desired statement.

\medskip
Next, we improve this sequential convergence to a full convergence as $t\to\infty$:

\begin{lemma}
Define the function $Q(t)$ by
\bea
Q(t)=\int_X(-R+\lambda+\tau)^2\o+\int_X |\na\tau|^2\o.
\eea
Then

{\rm (a)} There exists a sequence $t_j$, $j\leq t_j\leq j+1$, with $Q(t_j)\to 0$ as $j\to\infty$.

{\rm (b)} There exist  positive constants $A,\epsilon$ so that
\bea
{dQ\over dt}\leq -\epsilon \int_X ((\Delta\tau)^2+|\na R|^2)\o+ A Q(t).
\eea

{\rm (c)} We have $Q(t)\to 0$ as $t\to\infty$.
\end{lemma}

\medskip
\noindent
{\it Proof.} The statement (a) is a reformulation of the properties of the sequence $t_j$ that we just obtained, so we concentrate on the proof of (b). We evaluate
\bea
\p_t\int_X |\na\tau|^2\o
&=&\p_t\int_X \p_z\tau\p_{\bar z}\tau
=
2\int_X\p_z\dot\tau\p_{\bar z}\tau\nonumber\\
&=&
2\int_X \p_z(\kappa\Delta\tau-\tau(-R+\lambda+\tau))\p_{\bar z}\tau
\nonumber\\
&=&
-2\kappa\int_X (\Delta\tau)^2\o
+2\int_X \tau(-R+\lambda+\tau)\p_z\p_{\bar z}\tau.
\eea
Since $\tau$ is uniformly bounded, we can estimate the second term on the right hand side by
\bea
\left|\int_X \tau(-R+\lambda+\tau)\p_z\p_{\bar z}\tau\right|
&\leq& M\left|\int_X(-R+\lambda+\tau)\Delta\tau \right|\o
\nonumber\\
&\leq&
\epsilon \int_X|\Delta\tau|^2\o
+C_\epsilon \int_X(-R+\lambda+\tau)^2\o
\eea
and hence, for $\epsilon$ sufficiently small,
\bea
\p_t\int_X |\na\tau|^2\o
\leq -\kappa\int_X(\Delta\tau)^2\o+C \,Q(t).
\eea
Next,
\bea
\p_t\left(\int_X(-R+\lambda+\tau)^2\o\right)
&=&
\int_X 2(-\p_tR+\p_t\tau)(-R+\lambda+\tau)\o
+
\int_X (-R+\lambda+\tau)^2\p_t\o\nonumber
\eea
We estimate first the second term on the right hand side
\bea
\int_X(-R+\lambda+\tau)^2\p_t\o
=
\int_X(-R+\lambda+\tau)^3\o
\leq C\,Q
\eea
since $|R|$ and $\tau$ are bounded.
Next, the first term on the right hand side can be written as
\bea
&&2\int_X (\p_tR-\p_t\tau)(R-\lambda-\tau)\o \nonumber\\
&=&
2\int_X(\Delta R-(1+\kappa)\Delta\tau-(R-\tau)(-R+\lambda+\tau))(R-\lambda-\tau)\o \nonumber
\\
&=&
2\int_X(-|\na R|^2+\na R\cdot \na\tau)\omega-2(1+\kappa)\int_X\Delta\tau (R-\lambda-\tau)\omega
\nonumber\\
&+&
2\int_X(R-\tau)(R-\lambda-\tau)^2\o.
\eea
We can then write
\bea
2\int_X(-|\na R|^2+\na R\cdot\na\tau)\omega
\leq -\int_X|\na R|^2+\int_X|\na\tau|^2 \leq -\int_X|\na R|^2 + Q(t).
\eea
and
\bea
\left|2(1+\kappa)\int_X\Delta\tau (R-\lambda-\tau)\right|
&\leq& \epsilon \int_X(\Delta \tau)^2\o+C_\epsilon
\int_X(R-\lambda-\tau)^2\o\nonumber\\
&\leq&
\epsilon \int_X(\Delta\tau)^2\o+C_\epsilon Q(t).
\eea
Finally, since both $R$ and $\tau$ are uniformly bounded, we have
\bea
\left|2\int_X(R-\tau)(R-\lambda-\tau)^2\o\right|
\leq C\int_X \int_X(R-\lambda-\tau)^2\o\leq C\,Q(t).
\eea
The statement (b) is proved.

\smallskip
The statement (c) is an easy consequence of (b). Dropping the negative terms in the estimate for $dQ/dt$ and integrating from $t_j$ to any $t$ with $t_j\leq t\leq j+1$, we obtain
\bea
Q(t)\leq Q(t_j)e^{A(t-t_j)}
\leq e^A\,Q(t_j)
\eea
It follows that $Q(t)\to 0$ as $t\to\infty$, and (c) is proved.

\subsection{Convergence of the flow}

By Lemma \ref{prop:1} and Section \S  3.3, $\phi$ and $\tau$ are uniformly bounded. Thus we can apply Theorem 1 to get uniform $C^k$ bounds for $\phi$ and $\tau$ for any $k$.
It follows that for any sequence of times $t_j\to\infty$, there is a convergent subsequence in $C^\infty$
converging to $(\o_\infty,\tau_\infty)$. To get the full convergence instead of just subsequential convergence, it suffices to show that the limit is unique.

\medskip
This can be seen as follows. By the previous section, $\|\na\tau_j\|_{L^2(\o_0)}\to 0$, and hence $\tau_\infty$ is a constant. It follows that $\o_\infty$ is a metric of constant scalar curvature
\bea
Ric(\o_\infty)=(\lambda+\tau_\infty)\o_\infty
\eea
with volume 1. For $\chi < 0$, this specifies the metric $\o_\infty$ uniquely, and the convergence of the flow is established.

\subsection{A brief discussion of the case $\kappa=1$}

In this subsection we briefly discuss the $\kappa$-LYZ flow on compact Riemann surfaces with $\kappa=1$. Furthermore we assume that the initial data $\tau_0$ is positive as in Theorem \ref{RS}. Again by the maximum principle we know that $\tau$ remains positive along the flow, however we no longer have a quantity like in (\ref{eqn:r and tau}) to apply the maximum principle. The new idea here is to introduce an enhanced functional
$$\hat{E}(\phi,\tau)=E(\phi,\tau)+\int_X\frac{R^2}{\tau}\omega.$$
When $\kappa=1$, it is easy to see that
$$\frac{d}{dt}\hat{E}(\phi,\tau)=-\int_X(\dot\phi)^2\omega-\int_X\frac{|\nabla\tau|^2}{2\tau}\omega-\int_X\frac{2}{\tau}\left|\nabla R-\left(\frac{1}{2}+\frac{R}{\tau}\right)\nabla\tau\right|^2\omega\leq 0,$$
hence $\hat{E}(\phi,\tau)$ is monotone non-increasing along the flow. Given the initial data, by using a finer version of the Moser-Trudinger inequality, one can show that $\hat{E}(\phi,\tau)$ is uniformly bounded from below and the integral
$$\int_X\frac{R^2}{\tau}\omega$$
is uniformly bounded. Therefore if we assume that $\tau$ is uniformly bounded from above as in the second part of (\ref{eqn:assumption}), it follows that the Calabi energy $\int_XR^2\omega$ is uniformly bounded and we can apply Lemma \ref{struwe} to prove that $\phi$ is uniformly bounded. So we conclude that in the case of Riemann surfaces with $\kappa=1$ and $\tau_0>0$, the second assumption in (\ref{eqn:assumption}) implies the first, i.e., the evolving metric $\omega_t$ is uniformly equivalent to $\omega_0$.

\section{A Mabuchi-type Energy for the $\kappa$-LYZ Coupled Flow}
\setcounter{equation}{0}

In this section, we introduce a Mabuchi-type energy functional $M_{\omega_0}(\varphi,F)$ associated to the coupled system of scalar equations (\ref{eqn:varphi}) and (\ref{eqn:F}). If in addition to the assumptions in Theorem \ref{higherorder}, we assume that this functional $M_{\omega_0}$ is bounded from below along the flow, then we can establish the smooth convergence of (\ref{eqn:varphi}) and (\ref{eqn:F}), and hence of the metrics and closed forms in (\ref{LYZo}) and (\ref{LYZa-theta}). More precisely,

\begin{theorem}
\label{theorem3}
Under the same assumptions as Theorem \ref{higherorder} with $T=+\infty$ and $\lambda < 0$, if $M_{\omega_0}(\varphi_t,F_t)\ge -M_0$ for some $M_0>0$, then the solution $(\varphi_t, F_t)$ of the system (\ref{eqn:varphi}, \ref{eqn:F}) converges smoothly to a limit $(\varphi_\infty, F_\infty)$
satisfying the equations:
$$ \log\frac{\omega_{\varphi_\infty}^n}{\omega_0^n} = H_0 -\lambda \varphi_\infty + F_\infty,\quad \Delta_{\omega_{\varphi_\infty}} F_\infty = \tr_{\omega_{\varphi_\infty}} \alpha_0 - b.  $$
Here $\omega_{\varphi_\infty}=\o_0+i\ddb\varphi_\infty$ and $\alpha_\infty = \alpha_0  - i\ddb F_\infty$ are stationary solutions of the equations (\ref{LYZo}) and (\ref{LYZa-theta}), and $b$ is the constant defined in (\ref{b}).

Moreover, if $X $ does not admit any nontrivial holomorphic vector field, then the convergence is exponential.
\end{theorem}


\subsection{A Mabuchi-type energy functional}
Motivated by the Mabuchi $K$-energy, we introduce the following functional
for the coupled system (\ref{eqn:varphi}, \ref{eqn:F})
\footnote{In keeping with the traditional notation for the Mabuchi functional itself,
we have denoted the volume $\int_X \omega_0^n$ of $\o$ by $V$ instead of normalizing it to 1 as in the previous sections. This will also allow normalizing $\lambda$ to $-1$ in subsequent sections.}
\bea\label{eqn:K}
\mu_{\omega_0}(\varphi) = - \frac{n}{V}\int_0^1 \int_{X} \varphi_s' ( Ric(\omega_{\varphi_s}) -  \lambda \omega_{\varphi_s} - \alpha_0 ) \wedge \omega_{\varphi_s}^{n-1} ds,
\eea where $\varphi_s$ is a path in the space of K\"ahler potentials of $\omega_0$ connecting $0$ and $\varphi$, and $\omega_{\varphi_s} = \omega_0 + i\ddb \varphi_s$. It can be checked that the integral above is independent of the choice of the path $\varphi_s$. Furthermore, if $\varphi$ varies in a family of K\"ahler potentials, then we have
$$\frac{d}{dt}\mu_{\omega_0}(\varphi_t) = -\frac{n}{V}\int_X \dot\varphi ( Ric(\omega_{\varphi}) - \lambda \omega_{\varphi_t} - \alpha_0  )\wedge \omega_{\varphi_t}^{n-1}.$$
Let $\varphi$ and $F$ be the potentials of $\omega_t$ and $\alpha_t$ defined in (\ref{eqn:potential}). We define a functional associated with the flow equations (\ref{eqn:varphi}) and (\ref{eqn:F}) by
\bea\label{eqn:M functional}
M_{\omega_0}( \varphi, F  ) = \mu_{\omega_0}(\varphi)  - \frac{1}{V}\int_X F \omega_\varphi^n.
\eea
\begin{lemma}\label{lemma 10}
The functional $M_{\omega_0}$ defined in (\ref{eqn:M functional}) is monotone and non-increasing along the flow (\ref{eqn:varphi}) and (\ref{eqn:F}).
\end{lemma}
\noindent{\it Proof.}
This follows from direct calculations:
\bea\nonumber
\frac{d}{dt}M_{\omega_0}(\varphi, F) & = & -\frac{n}{V}\int_X \dot\varphi ( Ric(\omega_t) - \lambda \omega_t - \alpha_0) \wedge \omega_t^{n-1} - \frac{1}{V} \int_X (\dot F \omega_t^n + inF\ddb\dot\varphi \wedge \omega_t^{n-1})\\
& = & -\frac{n}{V}\int_X \dot\varphi ( Ric(\omega_t) - \lambda \omega_t - \alpha_0) \wedge \omega_t^{n-1} - \frac{n}{V} \int_X  \dot\varphi\, i\ddb F \wedge \omega_t^{n-1}\nonumber \\
& =\nonumber & \frac{n}{V}\int_X \dot\varphi\, i\ddb \dot\varphi \wedge \omega_t^{n-1}\\
& = & \nonumber  - \frac{1}{V} \int_X |\nabla \dot\varphi|^2_{\omega_t} \omega_t^n\le 0.
\eea

Next we recall the well-known $I$ and $J$ functionals of K\"ahler geometry,
$$I(\varphi) = \frac 1 V\int_X \varphi (\omega_0^n - \omega_\varphi^n)
=
{1\over V}\sum_{k=0}^{n-1}\int_Xi\p\varphi\wedge\bar\p\varphi\wedge\o_0^{n-1-k}\wedge\o_\varphi^k,
$$
and $$J(\varphi) = \frac{1}{V}\int_0^1 \int_X \varphi_s' (\omega_0^n - \omega_{\varphi_s}^n  )ds
=
{1\over V}\sum_{k=0}^{n-1}\int_X{n-k\over n+1}i\p\varphi\wedge\bar\p\varphi\wedge\o_0^{n-1-k}\wedge\o_\varphi^k,$$
where $\varphi_s$ is chosen as in the definition of $\mu_{\omega_0}(\varphi)$  in (\ref{eqn:K}). It is well-known that
$$I(\varphi)- J(\varphi)\ge \frac{1}{n+1}I(\varphi)\ge 0.$$
Recall that $H_0$ satisfies $Ric(\omega_0) - \lambda \omega_0 - \alpha_0 = i\ddb H_0$ and is normalized by $\int e^{H_0}\omega_0^n = V$. For any given path $\varphi_s$ connecting $\varphi$ and $0$, we define a family of smooth functions:
$$H_s = -\log\Big(\frac{\omega_{\varphi_s}^n}{\omega_0^n} \Big) - \lambda \varphi_s + H_0 + C(s) $$
where $C(s)$ is a normalizing function (constant in $X$ for each $s\in [0,1]$) such that $\int e^{H_s} \omega^n_{\varphi_s} = V$.  We can check that $H_s$ is a potential satisfying
$$Ric(\omega_{\varphi_s}) - \lambda \omega_{\varphi_s} - \alpha_0 = i\ddb H_s.$$ By the following calculation,
\bea\label{eqn:mu exp}
\mu_{\omega_0}(\varphi) &=& - \frac{n}{V}\int_0^1 \int_X \varphi_s'\, i\ddb H_s \wedge \omega_{\varphi_s}^{n-1} ds\\
& = & -\frac{n}{V}\int_0^1 \int_X H_s\, i\ddb \varphi'_s \wedge \o_{\varphi_s}^{n-1} ds\nonumber\\
& = & -\frac{1}{V}\int_0^1 \int_X H_s\frac{\partial}{\partial s} \omega_{\varphi_s}^n ds = -\frac{1}{V}\int_0^1\Big(\frac{d}{d s} \int_X H_s \omega_{\varphi_s}^n - \int_X \frac{\partial H_s}{\partial s} \omega_{\varphi_s}^n\Big) ds \nonumber \\
& = & \nonumber -\frac{1}{V}\int_X H_1 \omega_{\varphi}^n + \frac{1}{V}\int_X H_0 \omega_0^n + \frac 1V\int_0^1\int_X ( -\Delta_{\omega_{\varphi_s}} \varphi_s' -\lambda \varphi_s' + C'(s)    )\omega_{\varphi_s}^n ds\\
& = & \nonumber \frac 1 V \int_X \log\Big( \frac{\omega_\varphi^n}{\omega_0^n} \Big) \omega_\varphi^n + \frac{\lambda}{V} \int_X \varphi \omega_\varphi^n  - \frac 1 V \int_X H_0\omega_\varphi^n + \frac 1 V \int_X H_0 \omega_0^n -\frac{\lambda}{V}\int_0^1  \int_X \varphi_s' \omega_{\varphi_s}^n  ds\\
& = & \nonumber \frac 1 V \int_X \log\Big( \frac{\omega_\varphi^n}{\omega_0^n} \Big) \omega_\varphi^n - \lambda( I(\varphi) -  J(\varphi) )  - \frac 1 V \int_X H_0\omega_\varphi^n + \frac 1 V \int_X H_0 \omega_0^n,
\eea
we can see that $\mu_{\omega_0}(\varphi)$ is equal to the sum of an ``entropy part'', an ``energy part'' and some bounded terms.

\subsection{Convergence of the flow when $M(\varphi, F)>-\infty$ and $\lambda < 0$}
We can give now the proof of Theorem \ref{theorem3}.
Recall that we are making the basic assumption (\ref{eqn:assumption}) and also that $M(\varphi, F)\ge - M_0$ is bounded from below along the flow and $\lambda<0$. For notational simplicity we set $\lambda = -1$.

\begin{lemma}\label{lemma 11}
There exists a constant $C = C(\omega_0, \kappa, K)>0$ such that
$$-C \le \varphi_t + F_t \le C.$$
\end{lemma}
\noindent{\it Proof.}
We calculate
\bea\nonumber
\left(\frac{\partial }{\partial t} - \kappa \Delta \right)(F+\varphi  ) & = & -\kappa \tr_{\omega_t} \alpha_0 + \kappa b +  \log \frac{\omega_t^n}{\omega_0^n} - (\varphi + F) - \kappa n + \kappa \tr_{\omega_t} \omega_0,
\eea
the desired bound for $F + \varphi$ then follows from the maximum principle and the assumption (\ref{eqn:assumption}). The lemma is proved.

\medskip

Since $\dot \varphi + \varphi + F$ is uniformly bounded by the equation (\ref{eqn:varphi}) and the assumption (\ref{eqn:assumption}), we conclude from the above lemma that $\| \dot \varphi \|_{L^\infty}\le C(\omega_0, \kappa, K).$

\begin{lemma}\label{lemma 12}
There is a constant $C = C(\omega_0, \kappa, K, M_0) $ such that
$$\sup_X (|F_t| + |\varphi_t|)\le C,\quad \forall ~ t\in [0,T).$$
\end{lemma}

\noindent{\it Proof.}
From the expression of $\mu_{\omega_0} $ in (\ref{eqn:mu exp}) for $\lambda = -1$ we get $\mu_{\omega_0}(\varphi)\ge -C(H_0,\omega_0)$. Then the monotonicity of $M_{\omega_0}(\varphi, F)$ yields a uniform lower bound of $\int_X F_t \omega_t^n\ge -C(H_0,\omega_0,\varphi_0,F_0)$, where $\varphi_0$ and $F_0$ denote the initial values of $\varphi$ and $F$, respectively. On the other hand, since $I-J$ is given by linear combinations of integrals of the form $\int_X i\partial\varphi\wedge \bar \partial \varphi \wedge \omega_0^j \wedge \omega_\varphi^{n-1-j}$,
it is uniformly bounded by Lemma \ref{lemma 3}. The entropy part
\bea
\int_X \log( \frac{\omega_\varphi^n}{\omega_0^n}  )\omega_\varphi^n
\eea
is also bounded above by assumption. Thus the assumption that $M(\varphi, F)\ge - M_0$ implies a uniform upper bound of $\int_X F_t \omega_t^n \le C$. Combined with the gradient bound of $F$ in Lemma \ref{lemma 3}, this gives a $C^0$-bound for $F_t$, and by Lemma \ref{lemma 11}, we also get $\sup_X|\varphi_t|\le C$.
The lemma is proved.

\bigskip

Combining Lemma \ref{lemma 12} and Theorem \ref{higherorder}, we obtain a uniform estimate for all derivatives of $\varphi_t$ and $F_t$:
$$
\sup_X |\partial_t^k \nabla^l F_t   | + |\partial_t ^k \nabla^l \varphi_t  |\le C(k,l, \kappa, K,M_0),\quad\forall t\in [0, T).
$$
Recall that we are considering the case of infinite maximum existence time, $T = \infty$. From the proof of Lemma \ref{lemma 10}, we can deduce that
$$\int_0^\infty \int_X |\nabla\dot\varphi|^2 \omega_t^n dt <\infty,$$ and this implies that
$$
\lim_{i\to\infty} \int_i ^{i+2} \int_X |\nabla\dot\varphi|^2 \omega_t^n dt = 0.
$$
From the equation
$$(\frac{\partial}{\partial t} - \Delta)|\nabla \dot\varphi|^2= - |\nabla \nabla \dot\varphi|^2 - |\nabla\bar \nabla \dot\varphi|^2  - 2 Re( \< \nabla \dot\varphi, \bar \nabla( \dot\varphi + \dot F  )   \>   ) + |\nabla \dot\varphi|^2 - \alpha(\nabla\dot\varphi, \bar \nabla \dot\varphi),$$
we can use either Moser iteration or arguments similar to the ones in Section \S 3.4 to deduce that
$$\sup_{X\times[i+1,i+2]} |\nabla \dot\varphi|^2 \le C \int_{i}^{i+2} \int_X |\nabla\dot\varphi|^2 \omega_t^n dt\to 0,\quad {\mathrm{as}} ~ i\to\infty.$$ Therefore the smooth limit of $\dot\varphi_t$ as $t\to\infty$ must be constant. If $X$ does not admit any nontrivial holomorphic vector fields, we can use the argument of \cite{PSSW} to get the exponential convergence of $\int_X|\nabla\dot\varphi_t|^2\omega_t^n \le C e^{-\delta t}$ for some $\delta>0$ depending on the lower bound of positive eigenvalue of $- \Delta_{\omega_t}$ on vector fields in $T^{1,0}X$.

 By the $L^\infty$-bounds of $\varphi_t$ and $F_t$, we conclude that the smooth limits of $\dot\varphi_t$ and $\dot F_t$ must be zero. The limiting functions $\varphi_\infty$ and $F_\infty$ satisfy the equations
 $$\log \frac{\omega_{\varphi_\infty}^n}{\omega_0^n} = H_0 + \varphi_\infty + F_\infty, \quad \Delta_{\omega_{\varphi_\infty}} F_\infty - \tr_{\omega_{\varphi_\infty}} \alpha_0 + b = 0,$$
i.e. the $(1,1)$-forms $\omega_{\varphi_\infty} = \omega_0 + i\ddb \varphi_\infty$ and $\alpha_\infty = \alpha_0 - i\ddb F_\infty$ define stationary solutions of the flow (\ref{LYZo}) and (\ref{LYZa-theta}). The proof of Theorem \ref{theorem3} is complete.

\bigskip
Department of Mathematics, Columbia University, New York, NY 10027

tfei@math.columbia.edu, bguo@math.columbia.edu,
phong@math.columbia.edu



\end{document}